\documentclass[12pt, english]{article}
\usepackage[ruled]{algorithm2e}

\usepackage{cite}
\usepackage{amsmath, amssymb}
\usepackage{alltt}

\usepackage{subfig}
\usepackage{epsfig}

\usepackage{makeidx}         % allows index generation
\usepackage{graphicx}        % standard LaTeX graphics tool
\usepackage[T2A]{fontenc}
\usepackage[utf8]{inputenc}
\usepackage{babel}
\usepackage{url}
\newcommand{\vect} \mathbf

\begin{document}

\title{On fast matrix-vector multiplication with a Hankel matrix in  multiprecision arithmetics}
\author{Gleb Beliakov,
\\
School of Information Technology, Deakin University, \\221 Burwood Hwy, Burwood
3125, Australia \\ \texttt{gleb@deakin.edu.au}
}
%
%\date{}
\maketitle

\begin{abstract}
We present two fast algorithms for matrix-vector multiplication $y=Ax$, where $A$ is a Hankel matrix. The current asymptotically fastest method is based on the Fast Fourier Transform (FFT), however in multiprecision arithmetics with very high accuracy FFT method is actually slower than schoolbook multiplication for matrix sizes up to $n=8000$. One method presented is based on a decomposition of multiprecision numbers into sums, and applying standard or double precision FFT. The second method, inspired by Karatsuba multiplication, is based on recursively performing multiplications with matrices of half-size of the original. Its complexity in terms of the matrix size $n$ is $\Theta(n^{\log 3})$. Both methods are  applicable to Toeplitz matrices and to circulant matrices.
\end{abstract}


\baselineskip=\normalbaselineskip

%\thispagestyle{empty}

%as many elements of $\vect x$ smaller than or equal to  $x_i$ as those bigger than or equal to $x_i$, provided $n$ is odd.

\section{Introduction}

Matrix-vector multiplication $y=Ax$ is a fundamental step in many different algorithms, hence the interest in performing this operation as efficient as possible. One example is calculation of the eigenvalues of a matrix. In Lanczos algorithm \cite{Golub1996} matrix-vector multiplications are used iteratively to tridiagonalize a matrix, so that the subsequent application of QR iterations could compute the eigenvalues efficiently \cite{Bai2000}. Here we are interested in the case where the matrix $A$ is Hankel, although a similar approach covers the cases of Toeplitz  and circulant matrices  \cite{Heinig2011,Wang2006,Luk2000}.

For Hankel matrices the asymptotically fastest multiplication method is based on the Fast Fourier Transform (FFT), which has complexity $O(n \log n)$ as opposed to $O(n^2)$ for naive schoolbook multiplication. The  FFT is numerically stable \cite{Luk2000}, however its overheads mean that in terms of efficiency it offers advantages for sufficiently large matrix sizes $n$, which depend on the implementation and data types.

We are particularly interested  in performing operations in multiprecision arithmetics, with the accuracy of around ten thousand decimal places. The need to use such an accuracy came from the eigenvalues problem. It was reported \cite{Luk2000,Luk99analysisof} that Lanczos algorithm could  lead to a rather fast loss of numerical accuracy even for relatively small matrix sizes $n<64$. The loss of mutual orthogonality of the vectors generated in the Lanczos algorithm was blamed  as one of the main contributing factors of the loss of accuracy. In the motivating studies of Hankel matrix eigenvalues reported in \cite{YumatHankel}, computations with hundreds of decimal places using Mathematica package resulted in very significant losses of the accuracy (by a third of the initial accuracy) even for moderate matrix sizes up to $100$. Consequently our interest was to develop efficient approaches for matrix-vector multiplication in multiprecision arithmetics to serve as a key step of the Lanczos algorithm.

A straightforward approach is to use the FFT-based multiplication, but translate the FFT algorithm into multiprecision arithmetics. However it has several caveats. First, the computations in multiprecision arithmetics have significant costs which are disregarded when using standard single or double data type. For example, even a simple assignment operation involves the cost of copying a long array of bits representing the mantissa of a  multiprecision number. Furthermore, the multiplication operation has a significantly larger cost than addition, and hence the two operations are not equivalent in operations count. Addition has complexity $O(b)$ whereas the (asymptotically) fastest multiplication has complexity $O(b \log b)$, where $b$ is the number of bits in multiprecision numbers.
Second, multiprecision multiplication (for a large $b$) is itself based on FFT. Therefore FFT-based Hankel matrix multiplication would involve invocation of one FFT algorithm inside another, which may result in unforseen rises in the costs.

One approach we explore in this paper is to break down the multiprecision numbers into sums of double precision numbers and perform the standard FFT based multiplication of the decomposed matrix and vector. The advantage here is that standard highly optimised and often parallelised FFT routines can be used. The key ingredient in this method is an appropriate representation of the factors in the decomposed form, discussed later in this contribution.

An alternative fast matrix-vector multiplication method not based on FFT was recently discussed in \cite{Cariow}. This method based on Karatsuba multiplication \cite{Karatsuba} has computational complexity of $O(n^{\log 3})$ multiplications. However this approach involves a rather complicated housekeeping and involves a larger number of additions. Our method was also inspired by Karatsuba multiplication, however it is different from that of \cite{Cariow}, and is very simple to implement as a recursive process. Its complexity is also $O(n^{\log 3})$ in terms of multiplications and additions.

The rest of the paper is structured as follows.
in Section \ref{sec1} we present the necessary background. Section \ref{sec2} is dedicated to the description of the decomposition approach and FFT multiplication. Section \ref{sec3} describes Karatsuba-like multiplication and analyses its complexity. Section \ref{sec_conc} concludes.

\section{Preliminaries} \label{sec1}

A Hankel matrix has the following form
\begin{equation}
 A_H=
  \begin{bmatrix}
  a_1&a_2 &\dots & a_{n-1}& a_n\\
  a_2& a_3 & \dots& a_{n}&a_{n+1}\\
  \vdots & \vdots &\reflectbox{$\ddots$} &\vdots &\vdots\\
    a_{n-1}& a_{n}& \dots &a_{2n-3} &a_{2n-2}\\
    a_n& a_{n+1}& \dots &a_{2n-2} &a_{2n-1}
      \end{bmatrix}.
\label{AH}\end{equation}

For matrix size $n$ it can be represented through an array of elements of size $2n-1$, which we will also denote by $A$. A Toeplitz matrix has a different form
\begin{equation}
 A_T=
  \begin{bmatrix}
  a_n&a_{n+1} &\dots & a_{2n-2}& a_{2n-1}\\
  a_{n-1}& a_n & \dots& a_{2n-3}&a_{2n-2}\\
  \vdots & \vdots &\ddots &\vdots &\vdots\\
 a_{2}& a_3 & \dots& a_{n}&a_{n+1}\\
    a_1& a_{2}& \dots &a_{n-1} &a_{n}
      \end{bmatrix}.
\label{AT}\end{equation}
and can be converted to a Hankel matrix by reversing the order of rows or columns, i.e., by multiplying it by a permutation matrix. Both matrices can be embedded into circulant $2n \times 2n$ matrices which have the form
\begin{equation}
 A_C=
  \begin{bmatrix}
  a_1&a_n &\dots & a_{3}& a_2\\
  a_2& a_1 & \dots& a_{4}&a_{3}\\
  \vdots & \vdots &\ddots &\vdots &\vdots\\
    a_{n-1}& a_{n-2}& \dots &a_{1} &a_{n}\\
    a_n& a_{n-1}& \dots &a_{2} &a_{1}
      \end{bmatrix}.
\label{AT}\end{equation}

Fast multiplication of a circulant matrix by a vector can be achieved using FFT in $O(n \log n)$ operations \cite{Bai2000}. Consequently multiplication by a Hankel or Toeplitz matrix can also be performed in $O(n \log n)$ time using FFT.

The schoolbook multiplication of a Hankel matrix by a vector $y=A_Hx$ is computed as follows
$$
y_i = \sum_{i=1}^n a_{i+j-1} x_i.
$$

The FFT-based multiplication can be achieved by using direct ($FFT$) and inverse ($IFFT$) Fourier transforms
$$
\hat y=IFFT(FFT(\hat a) * FFT(\hat c)),
$$
where the auxiliary vectors are
$$\hat a=(a_n, a_{n+1},\ldots,a_{2n-1},a_1,\ldots a_{n-1})^T,
$$
$$
\hat x=(x_n, x_{n-1},\ldots,x_{1},0,\ldots, 0)^T,
$$
and
$$
\hat y=(y_1,y_2,\ldots,y_n,\ldots)^T,
$$
so that the desired product is given by the first $n$ components of $\hat y$.
The operation $*$ is the componentwise multiplication of two vectors. The operation count is $30 n \log n+O(n)$ additions and multiplications \cite{Luk2000,Luk99analysisof}. The authors of \cite{Luk2000,Luk99analysisof} claimed that the FFT algorithm becomes superior to the schoolbook multiplication for $n>16$ when using standard (hardware-based) data types. However for multiprecision data types this estimation does not hold. Our computations in multiprecision arithmetics show that for $b=32768$ bits accuracy, FFT becomes superior only for $n>8000$.

There are a  number of packages that offer arbitrary precision arithmetics, in particular Arprec, GMP, MPFR, Flint, Arb and some others \cite{GMPLIB,Fousse:2007:MMP,Arb,Arprec}. These packages implement the standard arithmetical operations and often a number of elementary and special functions.
However we were able to find only one implementation of FFT in multiprecision arithmetics, which is needed for the fast multiplication algorithm for Hankel matrices, namely the Multiprecision Computing toolbox for Matlab \cite{MatlabMPT}. This package is based on kissfft library \cite{kissfft}.

In multiprecision arithmetics one specifies the desired accuracy in terms of the number of bits $b$ used to hold the mantissa of the values. Multiprecision addition is relatively cheap with complexity $O(b)$, however multiplication is significantly more expensive. Depending on the number of bits used, one or another method is invoked. In particular for sufficiently large $b$ (depending on the hardware: for x386 architecture  the threshold for $b$ ranges between $2816$ and $7168$) FFT is used for multiplication in GMP and MPFR.

One approach to fast Hankel matrix multiplication is to change the basic type in a FFT library from float or double to multiprecision. This approach was taken in the Matlab Multiprecision computing toolbox \cite{MatlabMPT}. However there are also some alternatives which we wanted to explore. Since the fast multiprecision multiplication is itself based on FFT, it means that we need to perform one FFT as a basic step inside another. It might be possible to unroll these operations into just one standard FFT. On the other hand there are several parallel FFT libraries (in single or double precision) that make use of a) multiple cores of one CPU, b) many cores of a GPU, c) multiple CPUs using Message Passing interface (MPI). Changing the basic scalar type in these highly optimised libraries is not simple, as their performance heavily relies on the lengths of scalar type and machine words, extensive use of cache and other hardware related optimisations. The length of the scalar type may drastically worsen the performance.

\section{Decomposition approach}\label{sec2}

Bearing in mind that ideally we would want to make use of the existing optimised FFT libraries, we attempted to decompose each multiprecision number into a sum of float or double precision numbers, then applied a standard FFT library (including optimised parallel libraries), and eventually reconstructed multiprecision vectors from the arrays of floats. We call this approach decomposition. We now present the details of the decomposition approach to implementation of FFT in multiprecision arithmetics.

Every multiprecision number is represented as a sum $a=B_0 \sum\limits_{k=0}^l B^k a_{k}$, where $B$ is the base used, typically $B=2^{32}$ or $B=2^{64}$, and $B_0$ is a factor. The number $l$ is derived from the number of bits in the mantissa $b$ as $l=\lceil b/\log B \rceil$.

Let us show that we can write the product $y=A_H x$, in which every element of $A_H$ and $x$ is a sum $a_i=\sum_k a_{ik}$, $x_i=\sum_k x_{ik}$, as a product of an auxiliary larger Hankel matrix $\hat A$ by a vector $\hat x$, whose components are the terms of the sums. Consider the matrix
\begin{equation}
 \hat A=
  \begin{bmatrix}
  a_{11} &\dots &a_{1l} & a_{11} & \dots & a_{1l} & a_{21}  &\dots & a_{n,l-1}& a_{nl}\\
  a_{12} &\dots &a_{11} & a_{12}& \dots &a_{21}& a_{22} &\dots & a_{n,l}& a_{n+1,1}\\
  \vdots  &\vdots & \vdots &&\reflectbox{$\ddots$}&& \vdots & \vdots &\vdots &\vdots\\
    a_{nl}&  \dots  &\dots & \dots &\dots &\dots&\dots&\dots &a_{2n-1,l-1} &a_{2n-1,l}
      \end{bmatrix}.
\label{AH}\end{equation}
This matrix is a Hankel matrix determined by the vector $(a_{11}, \dots, a_{1l},  a_{11},  \dots , a_{1l},\dots ,a_{2n-1,l})$, and it has size $m=2nl$. Notice that the components of the sum for each element of $a_i$ appear twice. Let also
$$
\hat x=(x_{11},x_{12},\dots,x_{1l},0,\dots,0,x_{21},\dots,x_{2l},0,\dots,0,x_{31},\dots 0)^T.
$$
We can now compute the product $\hat y = \hat A \hat x$ by using FFT. Let us examine the element $\hat y_1$. It will be
$$
\hat y_1=a_{11}x_{11}+a_{12}x_{12} +\dots+a_{1l}x_{1l}+\dots = \sum_{k=1}^l \sum_{i=1}^n a_{ik}x_{ik}.
$$
Next
$$
\hat y_2 = \sum_{k=1}^l \sum_{i=1}^n a_{i,(k+1)mod\, l}x_{ik}, \;\hat y_3 = \sum_{k=1}^l \sum_{i=1}^n a_{i,(k+2)mod\, l}x_{ik},\;\dots
$$
$$
\hat y_l=a_{1l}x_{11}+a_{11}x_{12} +\dots+a_{1,l-1}x_{1l}+\dots = \sum_{k=1}^l \sum_{i=1}^n a_{i,(k+l-1)mod\, l}x_{ik}.
$$

Next we add the elements $\hat y_1 +\hat y_2 + \dots+\hat y_l$ and rearrange the sum to get
$$
\hat y_1  + \dots+\hat y_l = a_{11}\sum_{k=1}^l x_{1k} + a_{12}\sum_{k=1}^l x_{1k} +\dots+ a_{1l}\sum_{k=1}^l x_{1k}+a_{21}\sum_{k=1}^l x_{1k}+\dots=
$$
$$
\sum_{i=1}^n \sum_{k=1}^l \left( a_{ik} \sum_{j=1}^lx_{ij} \right)= \sum_{i=1}^n \left(\sum_{k=1}^l a_{ik}\right)\left( \sum_{j=1}^lx_{ij}\right)=\sum_{i=1}^n a_i x_i = y_1.
$$

Thus we obtained the first component of the desired solution $y$. When we proceed this way with the rest of the $l$-length segments of vector $\hat y$, we get the resulting formula
\begin{equation} \label{eq:decomp}
y_i = \sum_{k=l(i-1)+1}^{il} \hat y_k.
\end{equation}

Consequently we obtained the desired product $y=A_H x$ by decomposing each entry $a_i$ and $x_i$ into a sum, and then by multiplying the auxiliary matrix and vector $\hat y=\hat A \hat x$ constructed from the individual terms  of the sums, and recovering $y$ using partial sums of the elements of $\hat y$. Now, since every multiprecision number was represented as a sum of terms with standard accuracy, we can compute the product $A_H x$  by decomposing the entries into standard precision numbers, performing the product with standard precision, and then recovering the multiprecision result.

To analyse the complexity of this process, recall that the product $\hat A \hat x$ is performed using FFT, so the complexity will be $O(nl \log(nl))=O(nb \log(nb))$, depending on the number of bits $b$ of multiprecision numbers. This complexity is smaller than complexity of performing FFT in multiprecision arithmetics $O(nb \log n \log b))$, and the fact that we unrolled the two nested FFTs (we remind that multiplication in multiprecision arithmetics is also done by FFT for large $b$) into one FFT that uses standard data type may help to speed up the algorithm, since the existing FFT libraries are highly optimised. The wide availability of a number of such libraries on different platforms is a benefit.

On the other hand, a potential caveat is that losses of accuracy when performing computations with the standard accuracy could happen. The reason is that the terms of the sums representing mantissas of multiprecision numbers are very different in magnitude, and thus the standard FFT may lose accuracy during summations, even if the sum  \eqref{eq:decomp} is performed using higher accuracy. It remains to be seen whether the loss of accuracy occurs in practice, and to this end we plan to perform numerical experiments in a subsequent study.

\section{Iterative economic multiplication}\label{sec3}

In this section we present an alternative way of computing the product $y=A_H x$ inspired by Karatsuba multiplication \cite{Karatsuba}. Consider two rows of a Hankel matrix multiplied by vector $x$
\begin{eqnarray} \label{ex1}
&\quad\quad\quad\quad\quad\quad\quad\quad\ldots + a_i x_i + a_{i+1} x_{i+1} +\ldots \\ \nonumber
&\ldots + a_i x_{i-1} + a_{i+1} x_{i} +\ldots
\end{eqnarray}
Apparently we need four multiplications and two additions. However, consider a more economical scheme. Let
$$
C=(a_i + a_{i+1}) x_{i}, \; D=a_{i+1} (x_i - x_{i+1}), \; E= a_i (x_{i-1} - x_{i}).
$$
Then the first row in \eqref{ex1} can be written as $C-D$ and the second row as $C+E$. Hence we compute the same quantities using three multiplications and five additions. If we now perform these operations in every column of $A$ and every pair of rows, we achieve a saving of one quarter of multiplications, hence the operations count will be $\frac{3}{4} n^2$ multiplications, at the expense of a few (cheaper) additions. Furthermore the additions can be amortized, as exactly the same differences between $x_i$ and $x_{i+1}$ appear in the subsequent rows. We will perform a detailed operations count after we present the whole algorithm based on this observation.

Let us formalise the presented approach. We will create three auxiliary vectors (we remind that a Hankel matrix is represented by a $(2n-1)$ --vector). Let  $m=\lfloor \frac{n+1}{2} \rfloor$ and $m_1=\lceil \frac{n+1}{2} \rceil$. Let also
$$C: c_i=a_{2i-1}+a_{2i}$$ be a Hankel matrix containing the pairwise sums of neighbouring elements of $A$, $$D: d_i=a_{2i}, i=1,\ldots,m, \mbox{ and }E: e_i=a_{2i-1}, i=1,\ldots,m_1$$  be Hankel matrices containing the even and odd elements of (the vector representation of) $A$ respectively, and vectors  $$f: f_i=x_{2i-1}-x_{2i},\;\; g: g_i=x_{2i-2}-x_{2i-1}, \; i=1,\ldots,m_1$$ contain pairwise (even and odd) differences of the elements of $x$, with the convention $x_0=x_{n+1}=0$ and also $a_{2n}=0$, and vector $$h: h_i=x_{2i-1}, i=1,\ldots,m.$$

The Hankel matrices $C, D$ will be of size $m \times m$, and $E$ will be of size $m_1 \times m_1$ (for a odd $n$ we have $m=m_1$). Let us now compute the quantities
$p=C h$, $q=D f$ and $r=E g$. Now every odd numbered element of $y$ can be written as
$y_{2i-1}=p_i-q_i$ and every even numbered element $y_{2i}=p_i+r_i$, $i=1,\ldots,m$.

To see this, consider
$$
p_i = \sum_{j=1}^m c_{i+j-1} h_j = \sum_{j=1}^m (a_{i+2j-2} +a_{i+2j-1}) x_{2j-1}= \sum_{j=1}^m a_{(i-1)+2j-1} x_{2j-1} + \sum_{j=1}^m a_{(i-1)+2j} x_{2j-1},
$$
$$
q_i=\sum_{j=1}^m d_{i+j-1} f_{j} = \sum_{j=1}^m a_{(i-1)+2j} x_{2j-1} - \sum_{j=1}^m a_{(i-1)+2j} x_{2j},
$$
$$
r_i=\sum_{j=1}^{m_1} e_{i+j-1} g_{j} = \sum_{j=1}^{m_1} a_{(i-1)+2j-1} x_{2j-2} - \sum_{j=1}^{m_1} a_{(i-1)+2j-1} x_{2j-1}.
$$

We wee that $p_i-q_i$ collapses into
$$y_{2i-1}=\sum_{j=1}^m a_{(i-1)+2j-1} x_{2j-1} + \sum_{j=1}^m a_{(i-1)+2j} x_{2j}=\sum_{k=1}^n a_{(i-1)+k} x_{k},$$
and $p_i+r_i$ collapses into
$$y_{2i}=\sum_{j=1}^{m} a_{(i-1)+2j} x_{2j-1} + \sum_{j=1}^{m_1} a_{(i-1)+2j-1} x_{2j-2}=\sum_{k=1}^n a_{(i-1)+k+1} x_{k},$$
taking into account $x_0=0$ and $a_{2n}=0$.

Let us now count the number of operations. Vectors $p$ and $q$ require $m \times m$ multiplications and vector $r$ requires $m_1\times m_1$ multiplications, which is at most $m_1^2$. These vectors also require at most $m_1^2$ additions.

Construction of matrix $C$ requires $m$ additions and vectors $f$ and $g$ need  $m$ and $m_1$ additions respectively. So overall we need at most $3m_1^2$  multiplications and $3m_1^2+3m_1+n$ additions, which gives us at most $\frac{3}{4}(n+1)^2 + O(n)$ multiplications and additions, i.e., one quarter improvement over the schoolbook matrix-vector multiplication (and only $O(n)$  cost of  extra additions, which, as we see, are amortised).

Next, let us apply this algorithm recursively. We see that the products $p=C h$, $q=D f$ and $r=E g$ involve Hankel matrices of half the size of the original. The recursive algorithm is presented as Algorithm 1.

% Algorithm
\begin{algorithm}[!htp]
\SetAlgoNoLine
\KwIn{$n$, $a \in R^{2n-1}$, $x \in R^n$.}
\KwOut{$y\in R^n$.}
\begin{enumerate}
    \item[1] If $n = 2$ /* explicit  formula */ \\
        $y_1=a_1x_1+a_2x_2$; $y_2=a_2x_1+a_3x_2$; \\
        return $y$
  \item[2] /* Prepare auxiliary vectors */ \\
  $m=floor((n+1)/2)$, $m_1=ceil((n+1)/2)$ \\
  \item[3] For $i=1,\ldots,m$ do: $\;f_i=x_{2i-1}-x_{2i}$, $h_i=x_{2i-1}$. \hspace{12pt}/* take $x_{n+1}=0$ */
  \item[4] For $i=1,\ldots,2m-1$ do:$\;c_i=a_{2i-1}+a_{2i}$, $d_i=a_{2i}$. \hspace{5pt}/* take $a_{2n}=0$ */
  \item[5] For $i=1,\ldots,m_1$ do: $\;g_i=x_{2i-2}-x_{2i-1}$ \hspace{24pt} /* take $x_0=x_{n+1}=0$ */
  \item[6] For $i=1,\ldots,2m_1-1$ do: $\;e_i=a_{2i-1}$ \hspace{40pt}/* take $a_{2n+1}=0$ */
  \item[7] /* call Algorithm 1 recursively */ \\
  \begin{enumerate}
   \item[7.1] call Algorithm1($m,c,h, p$)
   \item[7.2] call Algorithm1($m,d,f, q$)
   \item[7.3] call Algorithm1($m1,e,g, r$)
    \end{enumerate}
   \item[8] For $i=1,\ldots,m$ do:
  \begin{enumerate}
  \item[8.1] $y_{2i-1}=p_i-q_i$
  \item[8.2] if $2i\leq n$ then $y_{2i}=p_i+r_i$
     \end{enumerate}
  \item[9] return $y$.
\end{enumerate}
\caption{Calculation of the product of a Hankel matrix by a vector. }
\label{alg:gauss}
\end{algorithm}

The analysis of complexity of the recursive algorithm is performed in a standard way using the recurrence relation $T(n)= 3 T(\frac{n}{2}) + f(n)$, where $T(n)$ is the cost of the Algorithm 1 for input of size $n$, $3$ is the number of invocations of this algorithm at Step 7, and $f(n)$ is the cost of preparation of the auxiliary vectors at Steps 3 -- 6, which is $\Theta(n)$. Then applying the master theorem \cite{Cormen2001} the cost of the recursive algorithm $T(n)=\Theta(n^{\log 3})\approx \Theta(n^{1.585})$, which is the same complexity as Karatsuba algorithm and the method reported in \cite{Cariow}.

The direct operations count can also be performed by noticing that when $n=2^k$, the depth of recursion is at most $k=\log n$. The algorithm stops recursion when the size of the subproblem is $2$, where there are two cases: 1) possibly $a_3 \neq 0$ when depth of recursion is $k-1$, and 2) $a_3=0$ when depth of recursion is $k$. In the latter case only three multiplications and one addition are needed. Then there will be at most $3 \times 3^{\lceil \log n\rceil}$ multiplications, which is no greater than $3 n^{\log 3}$. As far as additions are concerned, their number is no greater than
$3\lceil \frac{n+1}{2}\rceil (1+\frac{3}{2} + \frac{3}{2}^2 \ldots \frac{3}{2}^{k-1})+3^{\lceil \log n\rceil}+O(n)$. The first term, which is a partial sum of the geometric series giving $3\lceil \frac{n+1}{2}\rceil \frac{1-\frac{3}{2}^{k}}{1-\frac{3}{2}}\leq  6\lceil \frac{n+1}{2}\rceil\frac{3}{2}^{\log n}=6\lceil \frac{n+1}{2}\rceil n^{\log 3 -1}$ is the number of additions at Steps 3 -- 6, and the second term is the number of additions at Step 1. The total is no greater than
$  4n^{\log 3}+O(n)$ for $n>2$.

From the exact bound on the number of multiplication and additions it follows that the proposed Algorithm 1 will be computationally more efficient than the schoolbook matrix-vector multiplication for every $n>3$ regardless the cost of addition and multiplication as a function of the number of mantissa bits, and hence we expect this method to perform better in standard or multiprecision arithmetics.

The nature of the algorithm offers opportunities for parallelisation. Indeed the invocations of Algorithm 1 can be done by separate threads and on separate host computers. For example, a thread executing Algorithm 1 can spawn two other threads, so that the three recursive calls at Step 5 are executed in parallel. All threads will have the same amount of work and no load balancing issues should arise.

\section{Conclusion} \label{sec_conc}

We developed two algorithms for fast multiplication of a Hankel matrix by a vector, suitable in particular for multiprecsion arithmetics. The first algorithm decomposes the multiprecision numbers into sums of terms, and then uses fast Fourier transform to compute the product of the auxiliary Hankel matrix and vector of larger size. The complexity of this method is the same as FFT $O(n \log n)$, and the advantage is that existing FFT software that uses standard data types can be applied. However a potential problem is the loss of accuracy when using standard data types due to cancelation errors.

The second approach was inspired by Karatsuba multiplication and it is based on the reduction of the number of (multiprecision) multiplications from 4 to 3. When applied iteratively, the complexity of this approach becomes $\Theta(n^{\log 3})$, and the algorithm is very simple to implement. All computations are performed with full precision, so no additional loss of accuracy is expected. The algorithm is also easily parallelisable. The subsequent studies will evaluate the performance of both approaches empirically.

\bibliographystyle{plain}
\bibliography{median1}

\end{document}